\documentclass[12pt,a4paper]{article}

\usepackage[applemac]{inputenc}
\usepackage[english]{babel}
\usepackage[T1]{fontenc}

\usepackage{amsmath}
\usepackage{amsfonts}
\usepackage{amssymb}
\usepackage{amsthm}
\usepackage{appendix}
\usepackage{graphicx}
\usepackage{graphics}
\DeclareGraphicsExtensions{.jpg,.ps,.pdf,.png,.eps}
\usepackage{color}
\usepackage{booktabs}
\usepackage[pagebackref]{hyperref}
\usepackage{multirow}
\usepackage{url}
\usepackage{xcolor}

\usepackage {fancyhdr}
\usepackage{lastpage}
\usepackage{comment}

\usepackage{pdftricks}
\begin{psinputs}
   \usepackage{pstricks}
   \usepackage{multido}
\end{psinputs}

\usepackage[top=2.5cm, bottom=3cm, left=3cm , right=3cm]{geometry}

\usepackage{listings}

\newtheorem{thm}{Theorem}

\newtheorem{prop}[thm]{Proposition}
\newtheorem{lemme}[thm]{Lemma}
\newtheorem{rem}{Remark}

\newtheorem{definition}{Definition}

\newtheorem{hyp}{Assumption}
\newtheorem{corollaire}{Corollary}

\def\E{{\mathbb{E}}}
\def\R{{\mathbb{R}}}

\def\P{{\mathbb{P}}}

\def\N{{\mathbb{N}}}

\def\1{{\mathbf{1}}}
\def\F{{\mathcal{F}}}

\def\d{{\mathrm{d}}}

\def\ba{{\bar{a}}}

\def\la{{\langle}}
\def\ra{{\rangle}}
\def\bX{{\bar{X}}}
\def\bY{{\bar{Y}}}

\def\bN{{\bar{N}}}
\def\hN{{\hat{N}}}

\def\bmu{{\mu}}
\def\bm{{\bar{m}}}
\def\bC{{\bar{C}}}
\def\bD{{\bar{D}}}
\def\bA{{\bar{A}}}
\def\bJ{{\bar{J}}}
\def\bu{{\bar{u}}}
\def\bv{{\bar{v}}}

\def\tr{{\text{Tr}}}

\def\Var{{\mathrm{Var}}}

\def\da{\,{\partial_a}}
\def\ds{\,{\partial_s}}
\def\dt{\,{\partial_t}}

\def\dlambdak{\,{\partial_{y_k}}}

\def\bit{\begin{itemize}}
\def\eit{\end{itemize}}

\newcommand{\rmi}{{\rm (i) $\>\>$}}

\newcommand{\rmii}{{\rm (ii) $\hspace{1.5mm}$}}
\newcommand{\rmiii}{{\rm (iii)$\>\>$}}

\def\bc{\begin{center}}
\def\ec{\end{center}}

\def\bcom{}

\def\edoc{\end{document}}

\usepackage{natbib}
 \bibpunct[, ]{(}{)}{,}{a}{}{,}%
\newcommand{\cgblue}{\color{blue}}
\def\cred {\color{red}}

\title{Population viewpoint on Hawkes processes$^{1}$}
\author{Alexandre Boumezoued$^{2}$}
\begin{document}
\maketitle

\begin{center}
{\bf Abstract}
\end{center}
This paper focuses on a class of linear Hawkes processes with general immigrants. These are counting processes with shot noise intensity, including self-excited and externally excited patterns.  For such processes, we introduce the concept of age pyramid which evolves according to immigration and births. The virtue if this approach that combines an intensity process definition and a branching representation is that the population age pyramid keeps track of all past events. This is used to compute new distribution properties for a class of  linear Hawkes processes with general immigrants which generalize the popular exponential fertility function. The pathwise construction of the Hawkes process and its underlying population is also given.\\

{\bf Keywords:} Hawkes processes, branching, immigration, age pyramid, non-stationarity, laplace transform, thinning, Poisson point measure.

\footnotetext[1]{This work benefited from the financial support of the ANR project {\it Lolita} (ANR-13-BS01-0011),  and {\it Chaire Risques Financiers} of {\it Fondation du Risque}. }
\footnotetext[2]{Laboratoire de Probabilités et Modèles Aléatoires (LPMA), UMR CNRS 7599\\
Université Paris 6, 4 Place Jussieu, 75005 Paris, France.
Email: \href{alexandre.boumezoued@upmc.fr}{alexandre.boumezoued@upmc.fr}}

\section{Introduction}

This paper investigates the link between some population dynamics models and a class of Hawkes processes. We are interested in processes whose behavior is modified by past events, which are self-excited and externally excited. The introduction of a self-excited process with shot noise intensity is due to \cite{HAWKES1971} and the famous Hawkes process has been used until now for a variety of applications, including seismology, neuroscience, epidemiology, insurance and finance, to name but a few. The shot noise intensity of the Hawkes process $(N_t)$ is expressed as
$$
\lambda_t = \mu + \sum_{T_n<t} \phi(t-T_n),
$$ 
where the $T_n$ are the times of jump of the Hawkes process $N$ itself, $\mu> 0$ and $\phi$ is a non-negative function. In the Hawkes model, when an event occurs at time $T_n$, the intensity grows by an amount $\phi(t-T_n)$: this models the self-exciting property. Also, for many modeling purposes,  $\phi$ returns to zero as $t$ increases, so that the self-excitation vanishes after a long time. 
On the whole, each event excites the system as it increases its intensity, but this increase vanishes with time as it is natural to model the fact that very old events have a negligible impact on the current behavior of the process.
In the literature, more recent contributions focused on processes with self-exciting behavior and also some externally exciting component. To our knowledge, the Hawkes process with general immigrants has been introduced in \cite{BREMAUD2002}, and specific forms can also be found in recent studies motivated by financial applications, such as \cite{DASSIOS2011}, \cite{FILIMONOV2014} and \cite{RAMBALDI2014}, where external shocks, news arrivals and contagion are crucial to model.
In this paper, we are interested in a class of Hawkes processes with general immigrants (see \cite{BREMAUD2002}), whose intensity is of the form
$$
\lambda_t =\mu(t) + \sum_{T_n<t} \Phi_t(t-T_n,X_n) +  \sum_{S_k<t} \Psi_t(t-S_k,Y_k).
$$
In this model, the $T_n$ are the times of jump of $N$: if an event occurs for the system at time $T_n$, the intensity grows by an amount $\Phi_t(t-T_n,X_n)$, where $X_n$ is some mark. This part models the self-exciting property. In parallel, external events occur at times $S_k$ and excites the system of interest with some amount $\Psi_t(t-S_k,Y_k)$: this is the externally excited component.

Among the appealing properties of such models, one of them comes from the shot-noise form of the intensity. This is called the cluster (or branching) representation of the Hawkes process, and it is based on the following remark: if an event occured at time $T_n$, then $t-T_n$ is nothing but the "age" of this event at time $t$.  Few years later after the seminal work of \cite{HAWKES1971}, \cite{HAWKES1974} proposed the cluster representation of the self-exciting process. They interpreted it as an immigration-birth process with age: they proved that under some stationarity conditions, it can be described  as a branching Poisson process (also called Poisson cluster). Also, in \cite{DASSIOS2011}, a definition of a dynamic contagion process is given through its cluster representation.
Until now, most studies on the Hawkes process recalled the immigration-birth representation as follows: immigrants arrive at times given by a Poisson process with intensity $\mu$. Then each immigrant starts a new generation: it gives birth to new individuals with fertility function $\phi$, each one giving birth with same fertility function $\phi$. This is often used as a definition for the Hawkes process, providing a good intuition on its behavior.  The cluster representation of \cite{HAWKES1974} requires that the mean number of children per individual which is nothing but $\left\| \phi\right\|=\int_0^\infty \phi(a) \d a$ satisfies $\left\| \phi \right\|<1$. In our paper, we exhibit the underlying immigration-birth dynamics which does not require the stationary assumption. Each individual in the population has an age and a characteristic. The virtue if this approach that combines an intensity process definition and a branching representation is that the population age pyramid keeps track of all past events. This is used to compute new distribution properties for a class of  linear Hawkes process with general immigrants. 

In the literature, the distribution properties of the Hawkes process have first been studied under stationary conditions.  \cite{HAWKES1971} addressed second order stationary properties, whereas
\cite{ADAMOPOULOS1975} derived the probability generating functional under stationarity, by using the cluster representation of \cite{HAWKES1974}. In this work, \cite{ADAMOPOULOS1975} expressed the probability generating function as a solution to some functional equation. Furthermore, \cite{BREMAUD2002} introduced the framework for studying moments of the stationary Hawkes process by means of the Bartlett spectrum. Let us also mention two recent studies of the distribution properties under stationarity. The moment generating function has been expressed in \cite{SAICHEV2011} as a solution to some transcendental equation. 
In addition, \cite{JOVANOVIC2014} proposed a graphical way to derive closed form expressions for cumulant densities, leading to the moments of the stationary Hawkes process.  It is interesting to note that such recent contributions rely the stationary branching representation of \cite{HAWKES1974}.
Recently, the computation of statistical properties has gained attention under non-stationarity, both for mathematical analysis and statistical estimation techniques. However, the recent studies in this framework only focus on exponential fertility rates $\phi(t)=\alpha e^{\beta t}$. The tool they rely on is the infinitesimal generator of the intensity process $(\lambda_t)$ which is Markovian for such exponential fertility rate (see \cite{OAKES1975}). This includes the work of \cite{ERRAIS2010}, \cite{AIT2010}, \cite{DASSIOS2011},  and \cite{ZAATOUR2014}. Our paper generalizes these studies in a natural direction for a wider class of Hawkes processes.

\paragraph{Scope of this paper} The aim of this paper is
 \rmi to introduce the concept of age pyramid for general Hawkes processes and study its dynamics over time,
 \rmii to use this concept to compute new distribution properties for a class of fertility functions which generalize the popular exponential case, and 
 \rmiii to give a pathwise representation of the general Hawkes processes and its underlying immigration-birth dynamics. 
 We represent the population as a multi-type dynamics with ages, including immigration and births with mutations. Our population point of view that introduces the concept of age pyramid is inspired by \cite{workingpaper} (see also \cite{TRAN2008}). As highlighted in \cite{workingpaper}, the key idea is that the population structure in terms of ages and characteristics, which keeps track of past events, provides much more information than the intensity itself and allows to study the whole system. In this way, we address the computation of distribution properties of the Hawkes process with general immigrants for a wide class of time-dependent fertility functions. We also give the pathwise construction of the age pyramid represented as a measure-valued process solution to a stochastic equation driven by a Poisson point measure, which is the theoretical counterpart of the thinning numerical procedure.  Our approach seems to reconcile the two definitions of Hawkes processes, through an intensity process or a branching dynamics. 
  \\
The paper is organized as follows. 
Section \ref{section_population} focuses on the standard Hawkes process with time-independent fertility function. On this particular case, we give the population point of view and study the dynamics of the age pyramid over time. In Section \ref{section_application}, we use this concept to compute new distribution properties as moments and Laplace transform for a class of Hawkes processes which generalizes the popular exponential case.  Section \ref{section_construction} details the pathwise contruction of the standard Hawkes process and its underlying population. Our general  population representation and results are given in Section \ref{section_general_hawkes}, where we focus on Hawkes processes with general immigrants. In particular, we derive its dynamics and Laplace transform for a wide class of time-dependent fertility functions.

\section{Population point of view}
\label{section_population}

The definition of the (standard linear) Hawkes process through its intensity is given below.  Let $(\Omega, \mathcal{A}, \P)$ be a probability space satisfying the usual conditions. Recall that the intensity process $(\lambda_t)$ of a counting process $(N_t)$ is the $(\F_t^N)$-predictable process  such that $N_t- \int_0^t \lambda_s \d s$ is an $(\F_t^N)$- local martingale, where $(\F_t^N)$ denotes the canonical filtration of $(N_t)$.

\begin{definition}
\label{definition_hawkes} 
Let $\phi$ be a continuous and non-negative map. A Hawkes process $(N_t)$ with kernel $\phi$ is a counting process with canonical filtration $(\F_t^N)$ which admits an $(\F_t^N)$-predictable intensity
\begin{equation}
\label{equation_intensite}
\lambda_t= \bmu + \sum_{T_n<t} \phi(t-T_n) =  \bmu + \int_{(0,t)} \phi(t-s) \d N_s,
\end{equation}
where $\bmu>0$, and the $(T_n)$ are the times of jump of $(N_t)$.
\end{definition}

The previous definition provides the representation of the intensity process, which is interesting in order to study the behavior of the Hawkes process. But in fact, the whole information on the dynamics is lost. Indeed, it is interesting to go back to the branching representation of \cite{HAWKES1974} to have in mind the underlying population dynamics. First, immigrants arrive according to some Poisson process with parameter $\bmu$. Then each immigrant generates a cluster of descendants with the following rule: if an individual arrived or was born at some time $s$, it gives birth to new individuals with rate $\phi(t-s)$ at time $t$. Note that in fact, $t-s$ is nothing but the age at time $t$ of the individual born at time $s$. The birth mechanism can thus be reformulated as: any individual with age $a$ in the population gives birth with rate $\phi(a)$. The whole dynamics describes an immigration-birth process with age, in which the immigration rate is $\bmu$ and the birth rate is $\phi(a)$. 

Since the immigration-birth mechanism is crucial to understand the Hawkes dynamics, the aim now is to keep track of all ages in the population. One way to address this issue is to count the number of individuals with age below $\ba>0$ at time $t$, denoted $Z_t([0,\ba])$. This can be computed as the number of individuals arrived until time $t$ without those arrived before $t-\ba$, that is
$$
Z_t([0,\ba]) = N_t- N_{t-\ba}= \int_{(0,t]} \1_{t-s\leq \ba} \d N_s.
$$ 
The previous equation shows that with fixed $t$, this defines a measure on the space $\R_+$ of ages that is an image of the jump measure $\d N_t$. It can be written as
\begin{equation}
\label{equation_Zt}
Z_t(\d a)= \int_{(0,t]} \delta_{t-s}(\d a) \d N_s = \sum_{n=1}^{N_t} \delta_{t-T_n}(\d a). 
\end{equation}
Note that $Z_t(\d a)$ charges only $[0,t]$ since no individual born after time $0$ can reach an age greater than $t$. Formally, the measure $Z_t(\d a)$ puts a weight on the age of each individual alive at time $t$, therefore we call it age pyramid in reference to demographic analysis. In general, demographic studies focus on the number of individuals per age class of e.g. one year, so the quantity of interest is e.g. $Z_t([a,a+1))$.  The virtue of the measure representation is that one can compute a function $f$ of the population age structure by integrating it with respect to the age pyramid. To do this, we use the notation 
\begin{equation}
\label{equation_Ztf}
\la Z_t, f \ra = \int_{\R_+} f(a) Z_t( \d a)= \int_{(0,t]} f(t-s) \d N_s.%
\end{equation}
For example, the Hawkes process can be computed as $N_t=\la Z_t, \1 \ra$.
Also, the intensity process defined in Equation (\ref{equation_intensite}) can be rewritten using (\ref{equation_Ztf}) as
$$
\lambda_t=\bmu+ \la Z_{t-}, \phi \ra.
$$
The intensity is the sum of the migration intensity $\bmu$ and the individual birth intensities: this is indeed the intensity of an immigration-birth process with migration rate $\bmu$ and birth rate $\phi(a)$, in which all individuals behave independently.
Viewed as a stochastic process, $(Z_t(\d a))_{t\geq 0}$ is a measure-valued process. In fact, 
this age pyramid process, that is the measure-valued process $(Z_t(\d a))_{t\geq 0}$, is a Markov process (see \cite{THESETRAN}). Note however that its differentiation in time is not straightforward (see \cite{workingpaper} and Lemma \ref{lemme_semimartingale} below). The Markov property of the age pyramid process shows that all the information needed is contained in the population age structure. Let us mention the seminal point of view of \cite{HARRIS1963}, for who  "it does seem intuitively plausible that we obtain a Markov process, in an extended sense, if we describe the state of the population at time $t$ not simply by the number of objects present but by a list of the ages of all objects."  However, in practice this information is "too large" to perform tractable computations. In the next Section, we illustrate how to identify some minimal components to add to the Hawkes process  in order to make the dynamics Markovian. To do this, we first need to
address the time evolution of the age pyramid. The following Lemma details the dynamics of $\la Z_t,f \ra$ in the case where $f$ is differentiable. This is the key tool for our results in Section \ref{section_application}.

\begin{lemme}
\label{lemme_semimartingale}
For each differentiable $f: \R_+ \rightarrow \R$,
\begin{equation}
\label{equation_Ztff}
\la Z_t,f \ra =  f(0) \la Z_t, \1 \ra+ \int_0^t \la Z_s,  f' \ra \d s.
\end{equation}
\end{lemme}

\paragraph{Proof of Lemma \ref{lemme_semimartingale}.}
Let us  write between 
$s$ and $t$, $f(t-s)=f(0) + \int_s^t f'(u-s) \d u$ and use it into Equation (\ref{equation_Ztf}) to get
$
  \la Z_t,f \ra = f(0) \la Z_t, \1 \ra + \int_0^t \left( \int_s^t f'(u-s) \d u \right)\d N_s.
 $
 By Fubini's theorem, the last term of the sum is equal to 
$
\int_0^t \left( \int_0^u  f'(u-s) \d N_s \right) \d u, 
$
and by Equation (\ref{equation_Ztf}), this is equal to $\int_0^t \la Z_u,  f' \ra  \d u$. This concludes the proof. $\diamond$\\

The decomposition (\ref{equation_Ztff}) is classical in the field of measure-valued population dynamics (see \cite{TRAN2008} and \cite{workingpaper}). The first term refers to the pure jump part of arrivals of individuals with age $0$, whereas the second term of transport type illustrates the aging phenomenon: all ages are translated along the time axis. In particular, this shows why the intensity process $\lambda_t = \bmu + \la Z_{t-}, \phi \ra$ is Markovian in the case where the fertility function is exponential (see \cite{OAKES1975}), that is $\phi(a)=\alpha e^{\beta a}$. In this case, $\phi'= \beta \phi$, and Equation (\ref{equation_Ztff}) with $f \equiv \phi$ leads to the differential form
$$
\d \la Z_t,\phi \ra =  \alpha \d N_t+ \beta \la Z_t, \phi \ra \d t.
$$
Note that $\d N_t$ only depends on the past of $(\lambda_t)$ by means of the current value $\lambda_t$, which proves the Markov property.
This remark is the starting point of our study, which extends the exponential case in a natural setting.

\section{The exponential case generalized}
\label{section_application}

In this section, the aim is to use the concept of age pyramid process introduced in Section \ref{section_population} in order to compute several distribution properties for the non-stationary Hawkes process. In particular, we provide ordinary differential equations for first and second order moments and the Laplace transform.
All computations are performed under some assumption on the birth rate $\phi$ which naturally extends the popular exponential case.

\subsection{Assumption on the birth rate}

 \begin{hyp}
 \label{hypothese_equation_phi}
The map $a\in \R_+ \mapsto \phi(a)$ is non-negative, of class $\mathcal{C}^n(\R_+)$, and there exists $c=(c_{-1},...,c_{n-1}) \in \R^{n+1}$ such that $\phi$ statisfies
 \begin{equation}
 \label{equation_phi}
 \phi^{(n)}=c_{-1}+ \sum_{k=0}^{n-1} c_k \phi^{(k)},
 \end{equation}
 with  %
initial conditions $\phi^{(k)}(0)=m_k$, for $0 \leq k \leq n-1$. 
 \end{hyp}

The birth rates that satisfy Assumption \ref{hypothese_equation_phi} include the exponential case but also some fertility functions that are interesting for a variety of applications. 
Let us introduce the vector 
\begin{equation}
\label{equation_vecteurm}
m=(1,m_0, ..., m_{n-1})^T,
\end{equation}
 and the matrix $C=(C_{i,j})_{-1\leq i,j \leq n-1}$  given by $C_{i,i+1}=1$ for $0\leq i \leq n-2$ and $C_{n-1,j}=c_{j}$ for $-1\leq j \leq n-1$, all other components being zero. Since it is fully determined by the vector $c$, we denote
\begin{equation}
\label{equation_matriceC}
C(c)=
\begin{pmatrix}
0&  0& & &   \\
 & 0 & 1 &  &    \\
 &  & \ddots & \ddots &\\
 &  &  & 0 & 1\\
c_{-1} & c_0 & \cdots & c_{n-2} & c_{n-1} \\
\end{pmatrix}.
\end{equation}

 Equation (\ref{equation_phi}) can be rewritten $\Phi'=C\Phi$ where $\Phi=(1,\phi,...,\phi^{(n-1)})^T$, whose solution is given by $\Phi(a)=e^{aC}m$. Then $\phi$ can be recovered as the second component of the matrix $\Phi(a)$. In particular, if the polynomial $P(y)=y^n-\sum_{k=0}^{n-1} c_k y^k$ is split with distinct roots $y_1,...,y_p$ and corresponding multiplicity $n_1,...,n_p$, then $\phi$ can be written up to some constant as $\sum_{i=1}^p P_i(a) e^{y_i a}$ where $P_i$ is a polynomial with degree at most $n_i-1$. This is a sufficiently large set of functions to approximate any fertility function outside of the range of Assumption \ref{hypothese_equation_phi}.
As an example, the power law kernel is of importance for many applications. In the context of earthquakes, the Omori law describes the epidemic-type aftershock (ETAS) model: it corresponds to a specific form $\phi(a) \sim \frac{K}{a^{1+\epsilon}}$. Also in the field of financial microstructure, recent studies (see e.g. \cite{BOUCHAUD2013}) found that high-frequency financial activity is better described by a Hawkes process with power law kernel rather than exponential. %
The power law kernel with cut off can be approximated as in \cite{BOUCHAUD2013} up to a constant by the smooth function
$$
\phi(a)=  \sum_{i=0}^{M-1} \frac{e^{-a /(\tau_0 m^i)}}{(\tau_0 m^i)^{1+\epsilon}} - S e^{-a/(\tau_0 m^{-1})} ,
$$
where $S$ is such that $\phi(0)=0$. 
%
%
In general, one can use approximation theory to construct a sequence of fertility functions which tends to the original one. 
As a result, this constructs a sequence of Hawkes processes that approximate the original Hawkes process. 

\subsection{Dynamics}

Let us go back to the dynamics of the age pyramid over time. The key property that will allow us to compute distribution properties is that the population enables to identify the components to add to the Hawkes process and its intensity to make the dynamics Markovian. This is stated in the following proposition.

\begin{prop}
\label{proposition_markov}
Under Assumption \ref{hypothese_equation_phi}, the process $X_t= (\la Z_t, 1 \ra, \la Z_t, \phi \ra, ..., \la Z_t, \phi^{(n-1)} \ra )^T$ satisfies the dynamics
\begin{equation}
\label{equation_vectorielle}
X_t= N_t m + \int_0^t C X_s \d s,
\end{equation}
where we the vector $m$ and the matrix $C$ are given in (\ref{equation_vecteurm}) and (\ref{equation_matriceC}) respectively. %
In particular, $X$ is a Markov process.

\end{prop}

\paragraph{Proof of Proposition \ref{proposition_markov}.}

 %
Let us use Lemma \ref{lemme_semimartingale} to get for $0 \leq k \leq n-1$, with $f \equiv \phi^{(k)}$,
\begin{equation}
\label{equation_Ztphi_3}
\la Z_t,\phi^{(k)}\ra =    m_kN_t+ \int_0^t \la Z_s, \phi^{(k+1)} \ra \d s.
\end{equation}
By Assumption \ref{hypothese_equation_phi}, we get in particular
\begin{equation}
\label{equation_Ztphi_4}
\la Z_t,\phi^{(n-1)} \ra =   m_{n-1} N_t+ \sum_{k=-1}^{n-1} c_k \int_0^t \la Z_s, \phi^{(k)} \ra \d s,
\end{equation}
with convention $\phi^{(-1)} \equiv \1$.
This implies the dynamics (\ref{equation_vectorielle}) which also shows that %
$X$ is a Markov process. $\diamond$\\

The dynamics (\ref{equation_vectorielle}) for the $(n+1)$-dimensional vector $X_t$ gives a set of $n$ equations, the first component of $X$, which is the Hawkes process $N$, being free. In Section \ref{section_construction} we will give an equation on the Hawkes process $N$ by means of stochastic representation based on Poisson point measures. This will provide a full system of equations for the components of $X$ as well as a pathwise representation. For now, we are interested into deriving several distribution properties of the Hawkes process and its additional components in $X$.

\subsection{Moments}

\paragraph{First order moments}

The differential system of Equation (\ref{equation_vectorielle}) is linear and allows us to propose a straightforward differential equation for the first order moments. We also perform explicit computations for small dimensions $n=1$ and $n=2$.

 \begin{prop}
 \label{proposition_ordre1}
Under Assumption \ref{hypothese_equation_phi}, the vector map $u(t):=\E \left[ X_t \right]$ is solution to 
\begin{equation}
\label{equation_EDO_ordre1}
u'(t)=\bmu m + A u(t),
\end{equation} 
where 
the $(n+1) \times (n+1)$ matrix $A$ is given by 
\begin{equation}
\label{equation_matriceA}
A=C+mJ,
\end{equation}
where
\begin{equation}
\label{equation_vecteurJ}
J=(0,1,0,...,0),
\end{equation}
and the vector $m$ and the matrix $C$ are given in (\ref{equation_vecteurm}) and (\ref{equation_matriceC}) respectively.

\paragraph{Proof of Proposition \ref{proposition_ordre1}}
Let us use the martingale property of the compensated counting process, then use Fubini's theorem and the fact that Lebesgue measure charges no point to get
$
\E \left[N_t\right]=\int_0^t \left( \bmu + \E[\la Z_s, \phi \ra] \right) \d s. 
$
Now, let us take expectation in (\ref{equation_vectorielle}) and use the previous formula to get Equation (\ref{equation_EDO_ordre1}). 
$\diamond$

 \end{prop}

The differential equation (\ref{equation_EDO_ordre1}) allows to get explicit formulas for the expected number of events.  We recall the first order moment for the popular exponential case  $\phi(a)=e^{-ca}$ (see e.g. \cite{DASSIOS2011}) and also give the explicit formulas for the birth rate $\phi(a)=\alpha^2 a e^{-\beta a}$. Note that this case can be useful for a variety of applications to model a smooth delay at excitation. Remark also the different behavior of the first order moments, in particular in the critical case $\int_0^\infty \phi(a) \d a=1$, which corresponds to $c=1$ and $\alpha=\beta$.
For the two examples given below, the computations are left to the reader.

\begin{corollaire}
\label{corollaire1}
For the Hawkes process with $\phi(a)=e^{-c a}$, $c>0$, ($n=1$ in Assumption \ref{hypothese_equation_phi}), 
\begin{equation*}
\begin{split}
&\E[N_t]= \bmu \left(t + \frac{t^2}{2} \right) \; \text{ if } c=1,\\
&\E[N_t]=\frac{\bmu}{1-c} \left( \frac{e^{(1-c)t}-1}{1-c}-ct\right), \; \text{ if } c \neq 1.
\end{split}
\end{equation*}
\end{corollaire}

\begin{corollaire}
For the Hawkes process with  $\phi(a)=\alpha^2 a e^{-\beta a}$, $\alpha, \beta>0$, ($n=2$ in Assumption \ref{hypothese_equation_phi}),
\begin{equation*}
\begin{split}
&\E[N_t]=  \frac{\bmu}{8 \beta}\left(1-e^{-2\beta t} \right)+ \frac{3 \bmu}{4}t +\frac{\beta \bmu}{4}t^2, \;  \text{ if } \alpha=\beta,\\
& \E[N_t]=\frac{\bmu \beta^2}{\beta^2-\alpha^2}t + \frac{\alpha\bmu}{2} \left( \frac{e^{(\alpha-\beta)t}-1}{(\alpha-\beta)^2}  - \frac{e^{-(\alpha+\beta)t}-1}{(\alpha+\beta)^2} \right), \text{ if } \alpha \neq \beta.
\end{split}
\end{equation*}

\end{corollaire}

\paragraph{Second order moments}

In this subsection, we derive the dynamics of the variance-covariance matrix of the process $X_t:=(N_t, \la Z_t, \phi \ra, ..., \la Z_t, \phi^{(n-1)} \ra )^T$. As a consequence, we represent second order moments as the solutions to a linear ordinary differential equation. Our method is based on differential calculus with the finite variation process $(X_t)$ with dynamics (\ref{equation_vectorielle}) and could be extended to higher moments.

\begin{prop}
\label{proposition_second_order}
Let us introduce the variance-covariance matrix $V_t=X_t \bX_t$, where $\bX_t$ denotes the transpose of $X_t$. 
Then the matrix $V_t$ satisfies the dynamics
\begin{equation*}
\d V_t= \d N_t \left( X_{t-}\bm + m \bX_{t-} + m \bm \right)  + \d t \left( V_t \bC + C V_t \right).
\end{equation*}
In particular, the matrix $v(t)=\E\left[ V_t\right]$ satisfies the following ordinary differential equation
\begin{equation}
\label{equation_EDO_variance}
v'(t) =  v(t) \bA  + A v(t) + \bmu (m \bm + u(t) \bm + m \bar{u}(t)) +Ju(t) m \bm.
\end{equation}
where $u(t)$ is solution to (\ref{equation_EDO_ordre1})  and the matrix $A$ is defined in (\ref{equation_matriceA}).
\end{prop}

\paragraph{Proof of Proposition \ref{proposition_second_order}.}

Let us use the notation $X_t=(X^{[-1]}_t, X^{[0]}_t, ..., X^{[n-1]}_t )$.  Integration by parts leads to, for $-1 \leq l,k \leq n-1$,
\begin{equation*}
\d \left( X_t^{[k]} X_t^{[l]}\right)=X_{t-}^{[k]} \d X_{t}^{[l]} + X_{t-}^{[l]} \d X_t^{[k]} + m_k m_l \d N_t.
\end{equation*}
 The previous equation shows that
$
\d V_t=X_{t-} \d \bX_t  + (\d X_t) \bX_{t-} + \d N_t .m\bm.
$
By Proposition \ref{proposition_markov} and since Lebesgue measure charges no point we get
\begin{equation*}
\d V_t= \d N_t \left( X_{t-}\bm + m \bX_{t-} + m \bm \right)  + \d t \left( V_t \bC + C V_t \right).
\end{equation*}
Recall that $u(t)=\E\left[ X_t\right]$. Now, take expectation in the previous equation to get
\begin{equation*}
v'(t) =   \E\left[ (\bmu+ X_t^{[0]}) X_{t} \right]\bm + m \E \left[ (\bmu + X_t^{[0]}) \bX_{t} \right] + \left( \bmu +  \E [ X_t^{[0]}] \right) m \bm   +  v(t) \bC + C v(t).
\end{equation*}
Finally, note that $ X_t^{[0]}X_{t} = V_t \bJ$ where we recall that $J$ is defined by $J=(0,1,0,...,0)$, which makes the previous equation reduce to (\ref{equation_EDO_variance}). $\diamond$\\

We give explicit formulas for the popular exponential fertility function $\phi(a)=e^{-c a}$ and at a higher order for the case $\phi(a)=\beta^2ae^{-\beta a}$, which corresponds to the critical case since the mean number of children per individual satisfies $\int_0^\infty \phi (a) \d a=1$. Computations are based on the differential equation (\ref{equation_EDO_variance}) and are left to the reader.

\begin{corollaire}
\label{corollaire2}
For the Hawkes process with $\phi(a)=e^{-c a}$ ($n=1$ in Assumption \ref{hypothese_equation_phi}), 
\begin{equation*}
\begin{split}
&\Var(N_t)=\bmu t \left( 1+\frac{3}{2} t + \frac{2}{3} t^2 + \frac{1}{12} t^3\right) \; \mathrm{ if }\; c=1,\\
&\Var(N_t)=\frac{\bmu}{(1-c)^3} \left[ \frac{1-c/2}{1-c} e^{2(1-c)t} + \left( \frac{3c^2-1}{1-c} -2ct\right) e^{(1-c)t} - c^3 t + \frac{c(1/2-3c)}{1-c} \right], \; \text{ if } c \neq 1.
\end{split}
\end{equation*}
\end{corollaire}

\begin{corollaire}
 For the Hawkes process with $\phi(a)=\beta^2ae^{-\beta a}$, ($n=2$ in Assumption \ref{hypothese_equation_phi}), the variance of the intensity is given by
\begin{equation*}
\Var \left( \lambda_t \right) = \beta \mu \left(- \frac{7}{128} +\frac{3\beta}{32} t +\frac{\beta^2}{16} t^2 +\frac{1-\beta t}{8}e^{-2 \beta t}-\frac{9}{128}e^{-4\beta t}  \right).
\end{equation*}
\end{corollaire}

\subsection{Laplace transform}

The aim is to exhibit the exponential martingale associated with the process $X$ which consequently expresses its Laplace transform in a semi-explicit form. This is given in the following proposition. It is interesting to note that point \rmi refers to some forward martingale, whereas point \rmii focuses on backward martingality. Note that based on the Laplace transform, it is classical to recover moments of any order.

\begin{prop}
\label{proposition_laplace} Let us denote $(\F^X_t)$ the canonical filtration of the process $X$ and let us work under Assumption \ref{hypothese_equation_phi}.\\
\rmi For any deterministic and differentiable $A_t$, the following process is an $(\F^X_t)$-martingale:
\begin{equation}
\label{equation_exponentielle2}
\exp \left\{ A_t.X_t  - \int_0^t A_s.(CX_s) \d s- \int_0^t A_s'.X_s \d s -\int_0^t (e^{A_s . m} -1) \lambda_s \d s\right\}.
\end{equation}
\rmii For any $(n+1)$ real vector $v$,\\
\begin{equation}
\label{equation_laplace_prop}
\E\left[ \exp \left(v.X_T \right)\right] = \exp \left(-\bmu \int_0^T (1-e^{A_s.m})\d s \right),
\end{equation}
where the vector map $A$ satisfies the following non linear differential equation 
\begin{equation}
\label{equation_differentielle}
\bC A_t + A_t' +(e^{A_t.m}-1) J=0,
\end{equation}
with terminal condition $A_T=v$. Here, $v.X_T$ denotes the scalar product between $v$ and $X_T$, $J$ is defined in (\ref{equation_vecteurJ}), and $\bC$ is the transpose of the matrix $C$.\\
\rmiii Moreover, there exists a unique solution to Equation (\ref{equation_differentielle}).
\end{prop}

\paragraph{Proof of Proposition \ref{proposition_laplace}}
The exponential formula 
states that the following process is a martingale, for any deterministic $\alpha_s$,
\begin{equation}
\label{equation_exponentielle}
\exp \left\{ \int_0^t \alpha_s \d N_s -\int_0^t (e^{\alpha_s} -1) \lambda_s \d s\right\}.
\end{equation}
Now, by integration by parts and the use of Equation (\ref{equation_vectorielle}),
\begin{equation*}
\begin{split}
A_t.X_t &= \int_0^t A_s .\d X_s  + \int_0^t A_s'.X_s \d s\\ 
&= \int_0^t A_s. m \d N_s + \int_0^t A_s.(CX_s) \d s + \int_0^t A_s'.X_s \d s.
\end{split}
\end{equation*}
Then by Equation (\ref{equation_exponentielle}) with $\alpha_s=A_s.m$, the process in (\ref{equation_exponentielle2}) is a martingale.
To prove the second point, the aim is to find a martingale of the form $\exp \{ A_t.X_t +D(t)\}$ with some deterministic $D(t)$ and a terminal condition $A_T=v$. To do this, let us choose $A$ such that the random part in the integrant in (\ref{equation_exponentielle2}) vanish. Since $\lambda_s=\bmu + \la Z_{s-}, \phi \ra$, this amounts to get for each vector $X=( X^{[-1]},...., X^{[n-1]})$,
\begin{equation}
\bC A_t. X + A_t' .X+(e^{A_t.m}-1) (\bmu + X^{[0]}) =0.
\end{equation}
Let us now identify the term 
 in $X$, leading to the equation for $A$: $\bC A_t. X + A_t' .X+(e^{A_t.m}-1) X^{[0]}=0$, that is $\bC A_t + A_t' +(e^{A_t.m}-1) J=0$, where $J$ is defined in (\ref{equation_vecteurJ}).  If we set terminal condition $A_T=v$, we get Equation (\ref{equation_laplace_prop}). 
Finally, existence and uniqueness for Equation (\ref{equation_differentielle}) arises from Cauchy-Lipschitz theorem, since the map $Y\mapsto \bC Y+(e^{Y.m}-1) J$ is of class $\mathcal{C}^1$ on $\R^{n+1}$, and thus continuous and locally Lipschitz. $\diamond$\\

The previous result can be expressed in terms of a single function, and is derived below for the Hawkes process and its intensity. The proof is given in Appendix.

\begin{corollaire}
\label{corollaire_laplace}
Under Assumption \ref{hypothese_equation_phi}, the joint Laplace transform of the Hawkes process and its intensity is given for each real $\theta_1$ and $\theta_2$ by\\
\begin{equation}
\label{equation_laplace_prop}
\begin{split}
&\E\left[ \exp \left(\theta_1 N_T + \theta_2 \lambda_T \right)\right] = \exp \left\{-\bmu \left((-1)^{n} G^{(n)}(0) + \sum_{k=0}^{n-1} (-1)^{k+1} c_k G^{(k)}(0) \right)  \right\},
\end{split}
\end{equation}
where the function $G$ satisfies the non-linear ordinary differential equation: for each $0 \leq t \leq T$,
\begin{equation}
\label{equation_edo}
\begin{split}
 &(-1)^{n-1} G^{(n+1)}(t) + \sum_{k=0}^{n-1} (-1)^k c_k G^{(k+1)}(t)  + \exp \left( \theta_1 - c_{-1} G(t)+\sum_{k=0}^{n-1} b_k G^{(k+1)}(t)\right)-1 =0,\\
 & \text{ with terminal conditions } G^{(k)}(T)=0 \text{ for } 0 \leq k \leq n-1  \text{ and } G^{(n)}(T)=(-1)^{n-1}\theta_2,
\end{split}
\end{equation}
and for $0\leq k \leq n-1$, $b_k=(-1)^k \left( m_{n-1-k} - \sum_{l=k+1}^{n-1} m_{n-1-l} c_{n-l+k} \right)$.

\end{corollaire}


\section{Pathwise representation of Hawkes population}
\label{section_construction}

The aim of this Section is to detail the pathwise construction of the Hawkes process and its underlying population. This is done by means of stochastic differential equations driven by Poisson point measures. The virtue of this approach is that it seems to reconcile both the intensity process definition of the Hawkes process and its branching representation. We first describe the construction of the Hawkes process with reference to a Poisson measure, then exhibit the system of equations driving the generalized exponential case, and finally address the pathwise population dynamics for general birth rates.

\paragraph{Construction of the Hawkes process}

One question which arises with Definition \ref{definition_hawkes} refers to the construction of such process and the notion of pathwise uniqueness. An answer can be given by the thinning representation, %

which works as follows. Consider a Poisson point measure $Q(\d s, \d \theta)$ with intensity measure $q(\d s, \d \theta)=\d s \d \theta$ on $\R_+ \times \R_+$ (see e.g. \cite{CINLAR2011} for a definition), and denote $(\F_t^Q)$ the canonical filtration generated by $Q$. 
Note that the intensity measure $q$ is not finite and only $\sigma$-finite, which makes it impossible  to order the points in time of the Poisson point measure $Q$. But this flexible representation allows to represent a wide class of counting processes.
%
%
Let $(\lambda_t)$ be a $(\F_t^Q)$-predictable process such that a.s. for each $t>0$, $\int_0^t \lambda_s \d s <+\infty$. Then the following process $(N_t)$ is a counting process with $(\F_t^Q)$-predictable intensity $\lambda_t$:
$
N_t=  \int_{(0,t]} \int_{\R_+} \1_{[0,\lambda_s]}(\theta) Q(\d s, \d \theta).
$
Indeed, $N$ is clearly a counting process because each atom of $Q$ is weighted $1$ or $0$. Also, since a.s. $\int_0^t \lambda_s \d s <+\infty$, the martingale property for Poisson point measures ensures that 
$
N_t- \int_0^t \int_{\R_+} \1_{[0,\lambda_s]} ( \theta) \d \theta  \d s=N_t- \int_0^t \lambda_s \d s
$ is a $(\F_t^Q)$- local martingale. Now, let us describe the construction of the Hawkes process. Since the intensity in (\ref{equation_intensite}) is given as a particular form of the process itself, the idea is to define the Hawkes process as the solution to the stochastic equation 
\begin{equation}
\label{equation_thinning_hawkes}
N_t=  \int_{(0,t]} \int_{\R_+} \1_{[0, \bmu + \int_{(0,s)} \phi(s-u) \d N_u]}(\theta) Q(\d s, \d \theta).
\end{equation}
General results about existence and uniqueness for the Hawkes process (even non-linear) can be found in \cite{BREMAUD1996} and \cite{MASSOULIE1998} (see also \cite{FOURNIER2014} and the books of \cite{DALEY2008} and \cite{CINLAR2011}). The thinning method to represent a counting process as the solution of a stochastic equation is in fact classical. %
This general mathematical representation goes back to \cite{KERSTAN1964} and \cite{GRIGELIONIS1971}. One often refers to the thinning algorithms that have been proposed by \cite{LEWIS1978} and \cite{OGATA1981}, which are very useful to perform numerical simulations for quite complex intensity processes. A first advantage of the thinning formulation arises when one wants to show the existence of the Hawkes process. This is done by Picard iteration method (see \cite{MASSOULIE1998}): one constructs a sequence $(N^k)_{k \geq 0}$ of counting processes starting at $N^0\equiv 0$, and for $k \geq 0$,
\begin{equation}
\label{equation_cauchy}
N_t^{k+1}=  \int_{(0,t]} \int_{\R_+} \1_{[0, \bmu + \int_{(0,s)} \phi(s-u) \d N_u^k]}(\theta) Q(\d s, \d \theta).
\end{equation}
One can show that the sequence $(N^k)$ is Cauchy and thus converges to the desired process. Moreover, another advantage is to give strong uniqueness. With this issue, it appears that the thinning representation has the virtue to use "one noise once for all" and thus give pathwise construction and results. This is interesting to note  that this approach is used by \cite{FOURNIER2014} to show existence and uniqueness of an infinite graph of interacting Hawkes processes.
Due to the pathwise representation and the iterative construction, one can also identify each generation in the dynamics. Indeed, one sees in the construction of the Cauchy sequence in (\ref{equation_cauchy}) that $N^1$ counts the number of immigrants, whereas $N^2-N^1$ counts the children of immigrants, $N^3-N^2$ the grandchildren of immigrants, and so on. Generally, $N^{k+1}_t-N^k_t$ is the number of individuals in generation $k$ born before time $t$. This shows another advantage of the pathwise construction: what is called the "thinning parameter" $\theta$ gives additional information on the dynamics, making it possible in particular to study each generation separately. Before giving the representation of the age pyramid, we first go back to the extension of the exponential case.

\paragraph{Exponential case generalized}
We first address the particular case where the birth rate $\phi$ satisfies Assumption \ref{hypothese_equation_phi}. The dynamics of the $(n+1)-$dimensional vector $X_t:=(N_t, \la Z_t, \phi \ra, ..., \la Z_t, \phi^{(n-1)} \ra )^T$ is given in (\ref{equation_vectorielle}) by $\d X_t = \d N_t m + C X_t \d t$. This gives in fact $n$ equations, the first coordinate $N_t$ being free. The pathwise representation (\ref{equation_thinning_hawkes}) allows to derive the full system of equations by
$$
\d X_t = \int_{\R_{+}}  m \1_{\left[0, \bmu + X_{t-}^{[0]}\right]}(\theta) Q(\d t, \d  \theta) + C X_t \d t,
$$
where we recall the notation $X_t=(X^{[-1]}_t, X^{[0]}_t, ..., X^{[n-1]}_t )$.

\paragraph{Immigration-birth process with general fertility functions}

In the case where the birth rate is general, one has to represent the whole age pyramid, that is to give the thinning representation of the underlying immigration-birth process. In the field of population dynamics, this approach is used to construct extended birth-death processes with age in particular in \cite{workingpaper} (see also \cite{MELEARD2004} and \cite{TRAN2008}). 
From Equations (\ref{equation_Ztf}) and (\ref{equation_thinning_hawkes}), we get the pathwise representation
\begin{equation}
\label{equation_Zt_mesure}
Z_t( \d a) =  \int_{(0,t]} \int_{\R_+} \1_{[0, \bmu+ \la Z_{s-}, \phi \ra] } (\theta)  \delta_{(t-s)} ( \d a)Q(\d s, \d \theta).
\end{equation}
This illustrates the fact that the population at time $t$ is nothing but %
all individuals that arrived before time $t$ (immigration or birth); if an individual arrived at time $s$, its age at time $t$ is $t-s$. Note that in this form, the differentiation is not straightforward (see \cite{workingpaper}). But from Lemma \ref{lemme_semimartingale}, one can write the following (infinite) system of equations: for each differentiable $f: \R_+ \rightarrow \R$,
\begin{equation}
\label{equation_Ztf2}
\d \la Z_t,f \ra =  f(0) \int_{\R_+} \1_{[0, \bmu+ \la Z_{t-}, \phi \ra] } (\theta) Q(\d t, \d \theta) +  \la Z_t,  f' \ra \d t.
\end{equation}

This approach seems to reconcile the intensity process definition of the Hawkes process with its branching representation. Indeed, the population age pyramid is given through Equation (\ref{equation_Ztf2}) as a stochastic measure-valued process with its own intensity. Before going through the last Section on the Hawkes process with general immigrants, we briefly discuss the existing cluster representation in the following remark.

\begin{rem}%
We recall the definition of the Hawkes process in terms of a Poisson cluster introduced in \cite{HAWKES1974} and surveyed in the book of \cite{DALEY2003}.  Let $N_c(\d s)$ be a Poisson point measure on $\R_+$ with intensity measure $\bmu \d s$: this defines the cluster centers, also called ancestors. Let us introduce a family of point processes $\{\bN(\d t \mid s), s \in \R_+\}$. For each $s$, $\bN(\d t \mid s)$ defines the location of the offsprings within the cluster of an ancestor located at $s$. %
The cluster process $\hN$ counts the number of all offsprings of all immigrants by
$
\hN(\d t) = \int_{\R_+} \bN(\d t \mid s) N_c(\d s).
$
That is, the number of all offsprings up to time $t$ is given by 
$$
\hN([0,t])=\int_{\R_+} \bN([0,t] \mid s) N_c(\d s).$$
 Thus in the cluster representation, the Hawkes process can be written as the sum of the immigrants and their offsprings by
$$
N_c([0,t]) + \hN([0,t]).
$$

Note that the cluster representation has shown to facilitate the study of the Hawkes process under stationarity by using results in the field of branching processes. Our population representation seems to be the non-stationary counterpart, as it allows us to derive new distribution properties in this framework. 
Our population representation provides not only the size of the total progeny up to time $t$, but also a variety of quantities of interest depending on the population age structure. This has been used in Section \ref{section_application} in order to identify the components needed to make the dynamics Markovian. This will be also used in the following Section to study a class of  Hawkes process with general immigrants.

\end{rem}

\section{Towards more general Hawkes processes}%
\label{section_general_hawkes}

In this Section, we focus on a class of counting processes $N_t$ named as Hawkes processes  with general immigrants (see \cite{BREMAUD2002}), which is defined below.

\begin{definition}
\label{definition_hawkes_general}
A Hawkes process with general immigrants is a counting process $N_t$ whose intensity is given by%
\begin{equation}
\label{equation_intensite_hawkes_general}
\lambda_t =\mu(t) + \sum_{T_n<t} \Phi_t(t-T_n,X_n) +  \sum_{S_k<t} \Psi_t(t-S_k,Y_k),
\end{equation}
where the $T_n$ are the times of jump of $N$, the $S_k$ are the jumps of a counting process with deterministic intensity $\rho(t)$ and the $X_n$ (resp. $Y_k)$ are real positive iid with distribution $G$ (resp. $H$). The $(S_k)$, $(Y_k)$ and $(X_n)$ are assumed to be independent of each other.
\end{definition}

In this model, the $T_n$ are the times of jump of $N_t$: if an event occurs for the system at time $T_n$, the intensity grows by an amount $\Phi_t(t-T_n,X_n)$, where $X_n$ is some mark. This part models the self-exciting property. In parallel, external events occur at times $S_k$ and excites the system of interest with some amount $\Psi_t(t-S_k,Y_k)$: this is the externally excited component. The standard Hawkes process that have been studied in the previous Sections can be recovered by setting $\Phi_t(a,x)=\phi(a)$ and $\Psi_t(a,x)=0$.
The Hawkes process with general immigrants has been introduced and studied under stationary conditions by \cite{BREMAUD2002}. Due to their flexibility and natural interpretation, such models have gained recent attention for financial applications e.g. by  \cite{DASSIOS2011}, \cite{FILIMONOV2014} and \cite{RAMBALDI2014}. In particular, distribution properties of such process have been investigated by \cite{DASSIOS2011} in the case $\Phi_t(a,x)=\Psi_t(a,x)=xe^{-\delta a}$, in which framework the intensity process is Markovian. %
  The aim of this Section is to study the dynamics and characterize the distribution of the non-stationary Hawkes process with general immigrants for a larger class of fertility functions, possibly time-dependent, which extends the previous work of \cite{DASSIOS2011} in this direction. To do this, we represent a two-population immigration-birth dynamics with ages and characteristics.

\subsection{Description of the two-population dynamics}

The aim is to construct  populations of several individuals (or particles), each one having %
an age $a$ evolving over time, and a characteristic $x\in \R_+$. We construct two populations: the first one represents external shocks, whereas the second one represents events for the Hawkes process. \\
Each population $(i)$, $i=1$ or $2$, is represented at time $t$ as a measure which puts a weight on the age and characteristic of each individual, denoted $Z^{(i)}_t(\d a, \d x)$. The two populations are introduced based on Definition \ref{definition_hawkes_general} as
\begin{equation}
\label{equation_Zt1_mesure}
Z^{(1)}_t(\d a, \d x) = \sum_{S_k \leq t} \delta_{(t-S_k, Y_k)}(\d a, \d x) \text{ and } Z^{(2)}_t(\d a, \d x) = \sum_{T_n \leq t} \delta_{(t-T_n, X_n)}(\d a, \d x). 
\end{equation}
Since ages but also characteristics of individuals are involved, we prefer to call $Z^{(i)}_t$ population {\it structure} rather than {\it age pyramid}, which is more specific.
As for the standard Hawkes population representation, one can compute functions of the whole population structure, which can even depend on time. Consider a function $f_t(a,x)$ depending on time, and also on age and characteristics of individuals. This can be computed on the overall population by 
\begin{equation}
\label{equation_Zt2_mesure}
\la Z^{(i)}_t, f_t \ra=\int_{\R_+\times \R_+} f_t(a,x) Z^{(i)}_t(\d a, \d x),
\end{equation}
 for $i=1$ or $i=2$. For example, the Hawkes process is $N_t^{(2)}=\la Z_t^{(2)}, \1 \ra$.
Also, the intensity $\lambda_t$ of the Hawkes process $N_t^{(2)}$ given in Equation (\ref{equation_intensite_hawkes_general}) can be rewritten as
$$
\lambda_t= \mu(t)+ \la Z^{(2)}_{t-}, \Phi_{t} \ra +\la Z^{(1)}_{t-}, \Psi_{t} \ra.
$$
This shows that the underlying population dynamics works as follows.\\
\rmi Let us first describe the population $(1)$ of external shocks. It is made with immigrants that arrive in population $(1)$ with rate $\rho(t)$; at arrival, they have age $0$ and some characteristic $x$ drawn with distribution $H$. Any individual $(a,x)$ at time $t$ that belongs to population $(1)$ gives birth with rate $\Psi_t(a,x)$. The newborn belongs to population $(2)$; it has age $0$, and some characteristic drawn with distribution $G$.\\
\rmii Let us now complete the description of population $(2)$. In addition to births from population $(1)$, the population $(2)$ evolves according to two other kind of events: immigration and internal birth. Immigrants arrive in population $(2)$ with rate $\mu(t)$ with age $0$ and a characteristic drawn with distribution $G$. Any individual $(a,x)$ at time $t$ that belongs to population $(2)$ gives birth with rate $\Phi_t(a,x)$. The newborn also belongs to population $(2)$; it has age $0$, and some characteristic drawn with distribution $G$. 
\begin{figure}[h]
\centering
\label{image_contagion_3.pdf}
\includegraphics[scale=0.8]{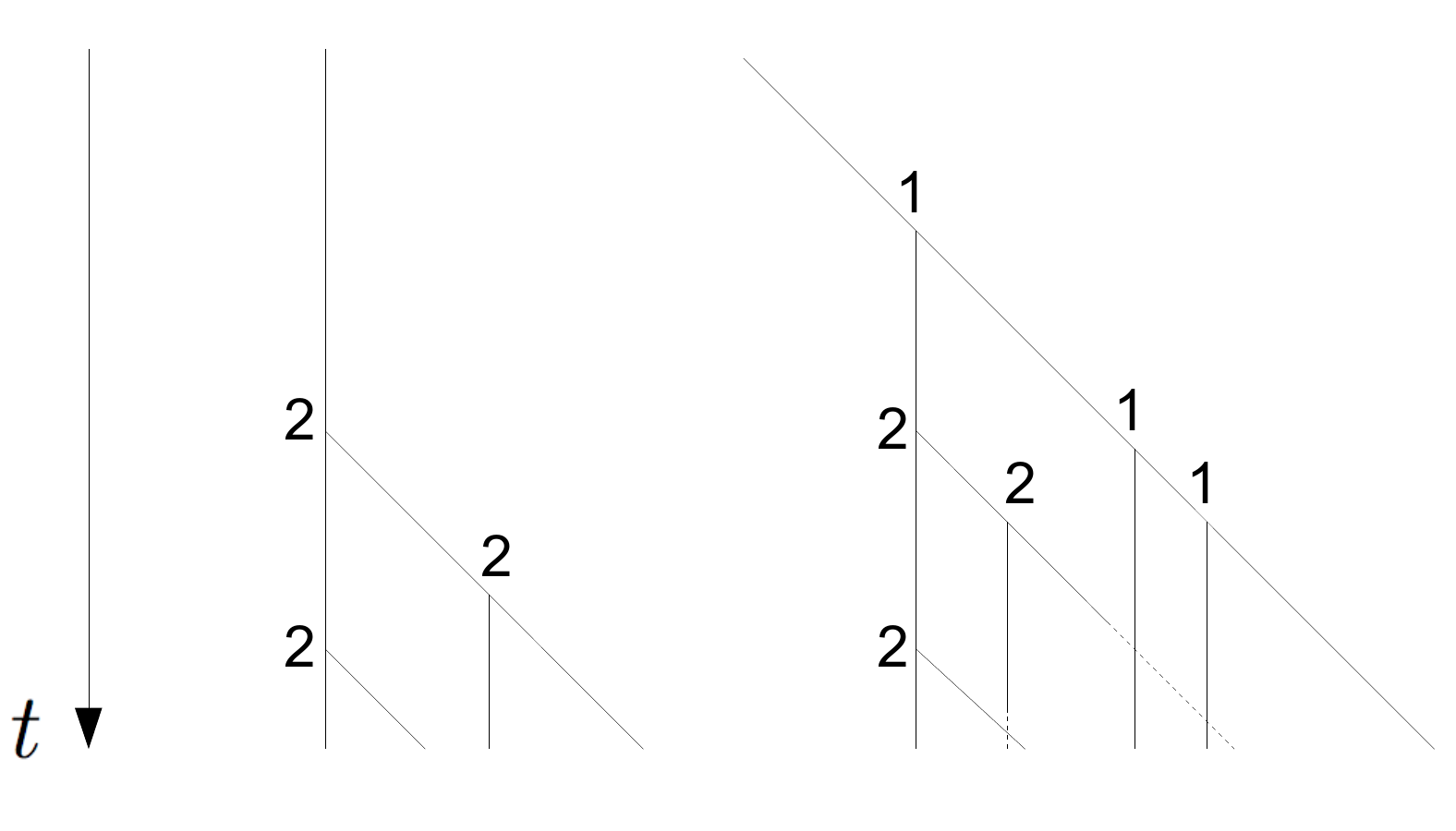}
\caption{Dynamics of the immigration-birth process: immigrants arrive in population $(1)$ (external shocks). Then each individual $1$ gives birth to individuals $2$ (events due to external shocks). In parallel, immigrants arrive in population $2$ (events due to the baseline intensity). Finally, each individual $2$ reproduce (self-excitation). The Hawkes process with general immigrants can be recovered as the number of individuals $2$. }
\end{figure}
This dynamics is illustrated in Figure 1.

As for our analysis of the standard Hawkes process, a crucial step is to study the dynamics of the population structure over time. That is, what is the dynamics of the process $\la Z_t^{(i)}, f_t \ra$ for $i=1$ or $2$ ? This is stated in the following lemma.

\begin{lemme}
\label{lemme_semimartingale2}
For each function $f: (t,x,a) \mapsto f_t(a,x)$ differentiable in $t$ and $a$, the dynamics of the process $\la Z_t^{(i)}, f_t \ra$ for $i=1$ or $2$ is given by
$$
\d \la Z_t^{(i)}, f_t\ra = \int_{\R_+} f_t(0,x) N^{(i)}(\d t, \d x) + \la Z_t^{(i)}, (\da + \dt) f_t \ra \d t,
$$
where the point measures $N^{(1)}$ and $N^{(2)}$ are given by 
\begin{equation}
\label{equation_N1N2}
N^{(1)}(\d t, \d x)=\sum_{k \geq 1} \delta_{(S_k,Y_k)}(\d t, \d x) \text{ and } N^{(2)}(\d t, \d x)=\sum_{n \geq 1} \delta_{(T_n,X_n)}(\d t, \d x).
\end{equation}

\end{lemme}

\paragraph{Proof of Lemma \ref{lemme_semimartingale2}}
The proof is a straightforward adaptation of that of Lemma \ref{lemme_semimartingale}, using (\ref{equation_Zt1_mesure}) and (\ref{equation_Zt2_mesure}) together with the fact that
$$
f_t(t-s)=f_s(0) + \int_s^t  (\da + \dt) f_u(u-s) \d u. \diamond
$$

Analogously to Lemma \ref{lemme_semimartingale}, this result exhibits the pure jump part in the left-hand side, whereas the drift part illustrates the aging term and the time-dependency.  The fact that the drift depends on both $\la Z_t^{(i)}, \da f_t\ra$ and  $\la Z_t^{(i)}, \dt f_t\ra$ is the starting point of our results derived in what follows.

\subsection{Main result}

In the following, we introduce the assumptions allowing to recover a finite dimensional Markovian dynamics.

\begin{hyp}
\label{hypothese_phi_psi}
\rmi The birth rates $\Phi$ and $\Psi$ are non-negative and satisfy $\Phi_t(a,x)=v(t)\phi(a,x)$ and $\Psi_t(a,x)=w(t)\psi(a,x)$, where
\begin{equation*}
 \phi^{(n)}(a,x)=c_{-1}+ \sum_{k=0}^{n-1} c_k \phi^{(k)}(a,x) \text{ and }  v^{(p)}(t)=d_{-1}(t)+ \sum_{l=0}^{p-1} d_l(t) v^{(l)}(t), 
 \end{equation*}
 with $n,p \geq 1$ and initial conditions $\phi^{(k)}(0,x)=\phi^{(k)}_0(x)$, and
 \begin{equation*}
 \psi^{(m)}(a,x)=r_{-1}+ \sum_{k=0}^{m-1} r_k \psi^{(k)}(a,x) \text{ and }  w^{(q)}(t)=k_{-1}(t)+ \sum_{l=0}^{q-1} k_l(t) w^{(l)}(t), 
 \end{equation*}
 with $m,q \geq 1$  and initial conditions $\psi^{(k)}(0,x)=\psi^{(k)}_0(x)$. 
 Note that we used the notation $f^{(k)}(a,x)= \da^{k} f(a,x)$.\\
\rmii The maps $(d_l)_{-1\leq l\leq p-1}$ and $(k_l)_{-1\leq l\leq q-1}$ are continuous.
\end{hyp}

\begin{rem} Assumption \ref{hypothese_phi_psi} defines a wide class of self and externally exciting fertility functions of the form $\Phi_t(a,x)=v(t) \phi(a,x)$. Let us first focus on the time-independent part and introduce $F(a,x)$ such that $F=(1,\phi,...,\phi^{(n-1)})^T$. Then $F'=C F$, where $C$ is defined in (\ref{equation_matriceC}). In particular, if the polynomial $P(y)=y^n-\sum_{k=0}^{n-1} c_k y^k$ is split with distinct roots $y_1,...,y_p$ and corresponding multiplicity $n_1,...,n_p$, then $\phi$ can be written up to some constant as $\sum_{i=1}^p P_i(x,a) e^{y_i a}$ where $P_i$ is a polynomial in $a$ with degree at most $n_i-1$ whose coefficients may depend on $x$. This includes the framework of \cite{DASSIOS2011} where $\Phi_t(a,x)=\Psi_t(a,x)=xe^{-\delta a}$.
As we also allow for time-dependency, such birth rates $\Phi$ and $\Psi$ that satisfy \ref{hypothese_phi_psi} seem also useful to define non-stationary Hawkes processes, and in particular to include seasonality. As an example, one can simply think of a kernel of the form $\cos^2(\alpha t) \phi(a,x)$ where $v(t)=\cos^2(\alpha t)$ satisfies $v''=4 \alpha^2 (1-v)$. 
\end{rem}

The aim of this part is to exhibit some exponential martingale which leads us to compute the Laplace transform of the whole dynamics.
This provides in particular the joint Laplace transform of the Hawkes process with general immigrants and its intensity. This is the main result of our paper. We first state the following Lemma.

\begin{lemme}
\label{lemme_matrice}
Let us define for $-1 \leq k \leq n-1$ and $-1 \leq l \leq p-1$, 
$X^{k,l}_t:=\la Z_t^{(2)}, \da^k \dt^l \Phi_t \ra$
and for $-1 \leq k \leq m-1$ and $-1 \leq l \leq q-1$,
$Y^{k,l}:=\la Z_t^{(1)}, \da^k \dt^l \Psi_t \ra.$
Let us also define the two matrices 
$$
M_t^{(2)}= \left( X^{(k,l)}_t \right)_{-1\leq k \leq n-1, -1 \leq l \leq p-1} \text{ and } M_t^{(1)}= \left( Y^{(k,l)}_t \right)_{-1\leq k \leq m-1, -1 \leq l \leq q-1}.
$$

\rmi Let us recall that $\bar{D}$ denotes the transpose of a given matrix $D$. The processes $M^{(1)}$ and $M^{(2)}$ follow the dynamics

\begin{equation}
\label{equation_matrice}
\d M_t^{(i)}= \int_{\R_+} W^{(i)} (t,x) N^{(i)} (\d t, \d x) + \left( C^{(i)} M_t^{(i)} + M_t^{(i)}  \bD_t^{(i)}  \right),
\end{equation}
where
\bit
\item $W_{k,l}^{(1)}(t,x)=w^{(l)}(t) \psi_0^{(k)}(x)$ for $-1\leq k \leq m-1$ and $-1\leq l \leq q-1$,
\item $W_{k,l}^{(2)}(t,x)=v^{(l)}(t) \phi_0^{(k)}(x)$ for $-1\leq k \leq n-1$ and $-1\leq l \leq p-1$,
\item $C^{(1)}=C(r)$, $C^{(2)}=C(c)$, $D^{(1)}_t=C(k(t))$ and $D^{(2)}_t=C(d(t))$ where $C(.)$ is defined by Equation  (\ref{equation_matriceC}).
\eit

\rmii As a consequence of the dynamics (\ref{equation_matrice}), $\left(M_t^{(1)},M_t^{(2)}\right)_{t \geq 0}$ is a Markov process.
\end{lemme}

\paragraph{Proof of Lemma \ref{lemme_matrice} }

We focus on the dynamics of the $X^{k,l}$, the problem being the same for the $Y^{k,l}$.
From Lemma \ref{lemme_semimartingale2}, for $0\leq k\leq n-2$ and $0\leq l\leq p-2$,
\begin{equation}
\label{equation_1}
\d X_t^{k,l}=v^{(l)}(t)\int_{\R_+} \phi_0^{(k)}(x) N^{(2)}(\d t, \d x) + (X_t^{k+1,l} + X_t^{k,l+1}) \d t.
\end{equation}

From Assumption \ref{hypothese_phi_psi}, $X^{n,l}_t = \sum_{k=-1}^{n-1} c_k X^{k,l}_t$ and $X^{k,n}_t = \sum_{l=-1}^{p-1} d_l(t) X^{k,l}_t$. This shows that for $0\leq l \leq p-2$,
\begin{equation}
\label{equation_2}
\d X_t^{n-1,l}=v^{(l)}(t)\int_{\R_+} \phi_0^{(n-1)}(x) N^{(2)}(\d t, \d x) + \left(\sum_{k=-1}^{n-1} c_k X^{k,l}_t + X_t^{n-1,l+1} \right) \d t,
\end{equation}
and for $0\leq k \leq n-2$,
\begin{equation}
\label{equation_3}
\d X_t^{k,p-1}=v^{(p-1)}(t)\int_{\R_+} \phi_0^{(k)}(x) N^{(2)}(\d t, \d x) + \left(X_t^{k+1,p-1} +  \sum_{l=-1}^{p-1} d_l(t) X^{k,l}_t \right) \d t.
\end{equation}
and also that
\begin{equation}
\label{equation_6}
\d X_t^{n-1,p-1}=v^{(p-1)}(t)\int_{\R_+} \phi_0^{(n-1)}(x) N^{(2)}(\d t, \d x) + \left(\sum_{k=-1}^{n-1} c_k X^{k,p-1}_t + \sum_{l=-1}^{p-1} d_l(t) X^{n-1,l}_t \right) \d t,
\end{equation}
In addition, Lemma \ref{lemme_semimartingale2} gives for $0\leq l \leq p-1$,
\begin{equation}
\label{equation_4}
\d X_t^{-1,l}=v^{(l)}(t) \d N^{(2)}_t + X_t^{-1,l+1} \d t.
\end{equation}
and for $0\leq k\leq n-2$,
\begin{equation}
\label{equation_5}
\d X_t^{k,-1}=\int_{\R_+} \phi_0^{(k)}(x) N^{(2)}(\d t, \d x) + X_t^{k+1,-1} \d t.
\end{equation}
Finally, by Assumption \ref{hypothese_phi_psi} again, we get the following two equations:
\begin{equation}
\label{equation_7}
\d X_t^{-1,p-1}=v^{(p-1)}(t) \d N^{(2)}_t + \left( \sum_{l=-1}^{p-1} d_l(t) X^{-1,l}_t \right)\d t,
\end{equation}
and
\begin{equation}
\label{equation_8}
\d X_t^{n-1,-1}=\int_{\R_+} \phi_0^{(n-1)}(x) N^{(2)}(\d t, \d x) + \left(\sum_{k=-1}^{n-1} c_k X^{k,-1}_t\right) \d t.
\end{equation}
From Equations (\ref{equation_1}) to (\ref{equation_8}), one then deduces the dynamics (\ref{equation_matrice}). $\diamond$

To ensure tractability of the Laplace transform derived in the following Theorem, we also state the following assumptions.
\begin{hyp}
\label{hypothese_existence}
For each $\lambda > 0$,
\begin{equation*}
\int_{\R_+} \exp \left( \lambda \max_{0 \leq k \leq n-1} \phi_0^{(k)}(x) \right) G(x) \d x < + \infty. 
\end{equation*}
\end{hyp}

Our main result is stated below. Note that the trace of the matrix $\bu M$ given by $\tr(\bu M)= \sum_{k,l} u_{k,l} M_{k,l} $  computes a linear combination of the components of a given matrix $M$, and recall that $\bu$ denotes the transposition of the matrix $u$.

\begin{thm}
\label{proposition_laplace_generale}
Let us denote $\mathcal{F}^M$ the canonical filtration generated by $(M^{(1)}, M^{(2)})$. Under Assumption \ref{hypothese_phi_psi}, \\
\rmi For any deterministic and differentiable matrix-valued $(A^{(1)}_t)$ and $(A^{(2)}_t)$ with derivatives  $(\mathring{A}^{(1)}_t)$ and $(\mathring{A}^{(2)}_t)$, the following process is an $\mathcal{F}^M$-martingale:
\begin{equation}
\label{equation_formule_exponentielle2}
\begin{split}
&\exp \bigg\{ \sum_{i=1}^2 \tr \left(  A_t^{(i)} M_t^{(i)} \right) %
- \int_0^t \tr \left(  A_s^{(i)} C^{(i)} M_s^{(i)} + A_s^{(i)} M_s^{(i)}  \bD_s^{(i)}  +  \mathring{A}_s^{(i)} M_s^{(i)} \right) \d s\\
&- \int_0^t \int_{\R_+} \left( e^{\tr \left(   A_s^{(1)} W^{(1)} (s,x) \right)}-1 \right) \rho(s) H(x) \d x \d s  \\
& - \int_0^t \int_{\R_+}  \left( e^{\tr \left(  A_s^{(2)} W^{(2)} (s,x) \right)}-1 \right) \left( \mu(s) + M_s^{(1)}[0,0] + M_s^{(2)}[0,0]  \right) G(x) \d x \d s \bigg\}.
\end{split}
\end{equation}

\rmii For each matrices $u$ and $v$ with dimensions $(n+1)(p+1)$ and $(m+1)(q+1)$ respectively, the joint Laplace transform can be expressed as
\begin{equation}
\label{equation_laplace_generale}
\begin{split}
\E \left[ \exp \left( \tr (\bu  M_t^{(1)} + \bv M_t^{(2)}) \right)\right]&=\exp \bigg\{ \int_0^t \int_{\R_+} \left( e^{\tr \left(   A_s^{(1)} W^{(1)} (s,x) \right)}-1 \right) \rho(s) H(x) \d x \d s  \\
& + \int_0^t \int_{\R_+}  \left( e^{\tr \left(  A_s^{(2)} W^{(2)} (s,x) \right)}-1 \right)  \mu(s)  G(x) \d x \d s \bigg\},
\end{split}
\end{equation}
where 
\begin{equation}
\label{equation_Ati}
\begin{split}
\text{ for } i\in \{1,2 \}, \;\mathring{A}_t^{(i)}  + A_t^{(i)} C^{(i)} + \bD_t^{(i)} A_t^{(i)} = \left\{  \int_{\R_+}  \left(1- e^{\tr \left(  A_t^{(2)} W^{(2)} (t,x) \right)} \right) G(x) \d x \right\} K,
\end{split}
\end{equation}

 with terminal conditions 
 \begin{equation}
\label{equation_conditions_terminales}
A_T^{(1)}=\bu \text{ and } A_T^{(2)}= \bv,
\end{equation}
where the matrix $K$ is given by $K=\bar{J} J$ and $J$ is given in (\ref{equation_vecteurJ}).
Moreover, solutions to  (\ref{equation_Ati})-(\ref{equation_conditions_terminales}) %
exist provided that Assumption $\ref{hypothese_existence}$ is satisfied.
\end{thm}

\paragraph{Proof of Theorem \ref{proposition_laplace_generale}}

 We begin by exhibiting the exponential martingale (\ref{equation_formule_exponentielle2}). Let us denote $\la N^{(i)}, H \ra_t = \int_0^t \int_{\R_+} H(s,x) N^{(i)}(\d s, \d x)$. For deterministic $\alpha(t,x)$ and $\beta(t,x)$, then by the classical exponential formula the following process is a martingale
\begin{equation}
\label{equation_formule_exponentielle}
\begin{split}
&\exp \bigg\{ \la N^{(1)} , \alpha \ra_t + \la N^{(2)}, \beta \ra_t - \int_0^t \int_{\R_+}  \left( e^{\alpha(s,x)}-1 \right) \rho(s) H(x)\d x \d s \\
& - \int_0^t \int_{\R_+}  \left( e^{\beta(s,x)}-1 \right) \left( \mu(s) + \la Z^{(1)}_{s-}, \Psi_{s} \ra + \la Z^{(2)}_{s-}, \Phi_{s} \ra \right) G(x) \d x \d s\bigg\}.
\end{split}
\end{equation}

The aim now is to compute the joint Laplace transform of the processes $M^{(1)}_t$ and $M^{(2)}_t$. This remains to compute $\E \left[ e^{\tr (\bu . M_t^{(1)} + \bv . M_t^{(2)})}\right]$, since  $\tr ( \bu . M) = \sum_{k,l} u_{k,l} M_{k,l}$.
Let us consider the two (deterministic) processes $A_t^{(1)}$ and $A_t^{(2)}$ with sizes $(m+1)(q+1)$ and $(n+1)(p+1)$ respectively. By integration by parts, $\d \left(  A_t^{(i)} M_t^{(i)} \right) =   A_t^{(i)} \d M_t^{(i)} + \mathring{A}_t^{(i)} M_t^{(i)} \d t$. From (\ref{equation_matrice}), we get the dynamics
\begin{equation*}
\begin{split}
&\d \tr \left( A_t^{(i)} M_t^{(i)} \right) =  \int_{\R_+} \tr \left(   A_t^{(i)} W^{(i)} (t,x) \right)N^{(i)} (\d t, \d x)\\
& + \tr \left(  C^{(i)} M_t^{(i)} + M_t^{(i)} \bD_t^{(i)}  +  \mathring{A}_t^{(i)} M_t^{(i)} \right) \d t
\end{split}
\end{equation*}
Let us now use Equation (\ref{equation_formule_exponentielle}) with $\alpha(t,x)= \tr \left(  A_t^{(1)} W^{(1)} (t,x) \right)$ and  $\beta(t,x)= \tr \left(   A_t^{(2)} W^{(2)} (t,x) \right)$ to get the martingale $(\ref{equation_formule_exponentielle2})$.

To get the Laplace transform, it remains to make the random part of the integrant in (\ref{equation_formule_exponentielle2}) vanish. To do this, let us first identify the term in $M^{(1)}$ to get the linear equation (\ref{equation_Ati}) for $i=1$.
In addition, the term in $M^{(2)}$ leads to (\ref{equation_Ati}) for $i=2$. 
If we set terminal conditions (\ref{equation_conditions_terminales}), we get the Laplace transform (\ref{equation_laplace_generale}) by the martingale property of (\ref{equation_formule_exponentielle2}).

To conclude on the existence and uniqueness, we use Cauchy-Lipschitz theorem. To show that solution of class $\mathcal{C}^1$ to  (\ref{equation_Ati}) %
exist and is unique, it is sufficient to prove that the map $(Y,t) \mapsto \int_{\R_+}  e^{\tr \left( Y W^{(2)} (t,x) \right)}G(x) \d x$ is of class $\mathcal{C}^1$. Since the integrant is $\mathcal{C}^1$ by Assumption \ref{hypothese_phi_psi} \rmi and \rmii, it is sufficient to prove that its gradient given by 
\begin{equation}%
\label{equation_gradient}
\left( e^{\tr \left( Y W^{(2)} (t,x) \right)} \bY, e^{\tr \left( Y W^{(2)} (t,x) \right)} \partial_t W^{(2)} (t,x) \right)
\end{equation}
is locally bounded by some quantity that is independent of $Y$ and $t$, and is integrable with respect to $G$.
Let us use some localization argument, and define the set $B(0,r)=\{A  \text{ real } (n+1)\times(p+1)\text{ matrix such that}  \;\left\| A\right\|_\infty \leq r \}$, where $r>0$ and $\left\| A\right\|_\infty = \max_{-1 \leq i \leq n-1} \sum_{j=-1}^{p-1} \left| A_{i,j} \right| $. Now, for $(Y,t) \in B(0,r)\times [0,T]$ we get
\begin{equation*}
\begin{split}
&\exp \left( \tr \left( Y W^{(2)} (t,x) \right) \right) \\
&\leq \exp \left(  \sum_{i=-1}^{n-1} \sum_{k=-1}^{p-1} \left| Y_{i,k} \right| \left| W^{(2)}_{k,i} (t,x)\right|\right)\\
&\leq \exp \left( (n+1) \max_{-1\leq i \leq n-1} \sum_{k=-1}^{p-1}  \left| Y_{i,k} \right| \left| W^{(2)}_{k,i}  (t,x)\right|  \right)\\
& \leq \exp \left( r (n+1) \max_{-1\leq l \leq p-1} \sup_{t\in [0,T]} \left| v^{(l)}(t) \right| \max_{-1\leq k \leq n-1} \left| \phi_0^{(k)}(x) \right|  \right),
\end{split}
\end{equation*}
where the last inequality uses that $Y\in B(0,r)$.
As for the first component of (\ref{equation_gradient}), $\left| \bY_{l,k} \right| \leq r$ and for the second component we have \\
$\left| \partial_t W^{(2)}_{k,l}(t,x)\right| \leq \left| \phi_0^{(k)}(x) \right| \sup_{t \in [0,T]} \left| v^{(l+1)}(t) \right|$, this concludes the proof by the use of Assumptions \ref{hypothese_phi_psi} and \ref{hypothese_existence}.
$\diamond$

\subsection{On the pathwise representation}

As for the standard Hawkes process, it is possible to give a pathwise representation of the Hawkes process with general immigrants and its underlying population. To do this, let us first extend the thinning construction in Section \ref{section_construction} to point processes with marks. Poisson point measures can be used, not only to represent counting processes, but also general random point measures on $\R_+ \times E$, say $\Gamma(\d s, \d y)= \sum_{n\geq 1}\delta_{(T_n,Y_n)} (\d s,\d y)$, where $(E, \mathcal{E})$ is some measurable space. As for the Hawkes process, $T_n$ is seen as the time at which an individual arrives in the population (immigration or birth). In addition, $E$ represents the space of characteristics and the mark $Y_n$ refers to the characteristic inherited by the individual that arrived at time $T_n$. Let us  construct a random point measure $\Gamma(\d s, \d y)$ with general intensity measure $\gamma(\d s, \d y)$ assuming that it admits a density: $\gamma(\d s, \d y)= \gamma(s,y)  \, \d s \, \mu(\d y)$.  In this model, events occur with intensity $s \mapsto \int_{x \in E} \gamma(s,x) \mu(\d x)$, and if a birth occurs at time $T_n$, then the characteristics $Y_n$ of the newborn are drawn with distribution $\frac{\gamma(T_n, y) \mu(\d y) }{\int_{x \in E}\gamma(T_n, x) \mu(\d x)}  $.\\
Let $Q(\d s, \d y, \d \theta)$ be a Poisson point measure on $\R_+ \times E \times \R_+$ with intensity measure $\d s  \mu( \d y)\d \theta$. Let us still denote $(\F^Q_t)$ the canonical filtration generated by $Q$, and introduce $P(\F^Q_t)$ the predictable $\sigma$-field associated with $\F^Q_t$. We further assume that $\gamma(t,y)$ is $P(\F^Q_t) \times \mathcal{E}$-measurable and also that $\int_0^t \int_E \gamma(s,y) \d s\mu( \d y) <+\infty$ a.s.. Now, define 
\begin{equation}
\label{equation_representation}
\Gamma(\d s, \d y)= \int_{\R_+}  \1_{[0, \gamma(s,y)]}(\theta) Q(\d s, \d y, \d \theta).
\end{equation}
This clearly defines a point measure and the martingale property for $Q$ ensures that the random point measure $\Gamma(\d s, \d y)$ has intensity measure $\gamma(s,y) \d s \mu( \d y)$.
Such construction can be found in \cite{MASSOULIE1998}; we refer to this paper for more details. 

We are now ready to construct the two point measures $N^{(1)}$ and $N^{(2)}$ given in Equation (\ref{equation_N1N2}).  
Let us introduce two independent Poisson point measures $Q^{(1)}(\d t, \d x, \d \theta)$ and $Q^{(2)}(\d t, \d x, \d \theta)$ on the probability space $(\Omega,\F, \P)$ (enlarged if necessary) with same intensity measure $\d s \d x \d \theta$ on $\R_+ \times \R_+ \times \R_+$. The first point measure is immediate to construct since its intensity does not depend on it. Indeed, one can define
\begin{equation*}
N^{(1)}( \d t, \d x) =  \int_{\R_+} \1_{[0, \rho(t) H(x)]} (\theta) Q^{(1)}(\d t, \d x, \d \theta).
\end{equation*}
We emphasize that this is not an equation on $N^{(1)}$ since its intensity does not depend on $N^{(1)}$ itself.
As for the second point process related to the Hawkes process, the intensity is given as a particular for of the process itself. Indeed, the intensity measure of the point measure $N^{(2)}(\d t, \d x)$ is given by $\lambda_t G(x)$ where $\lambda_t$ can be written using (\ref{equation_intensite_hawkes_general}) as
$$
\lambda_t= \mu(t) + \int_{(0,t)} \Phi_t(t-s,x) N^{(2)}(\d s, \d x) +  \int_{(0,t)} \Psi_t(t-s,x) N^{(1)}(\d s, \d x).
$$
Then the point measure $N^{(2)}$ can be defined as the solution to the following equation:
\begin{equation*}
N^{(2)}( \d t, \d x) =  \int_{\R_+} \1_{\left[ 0,  \left(\mu(t) + \int_{(0,t)} \Phi_t(t-s,x) N^{(2)}(\d s, \d x) +  \int_{(0,t)} \Psi_t(t-s,x) N^{(1)}(\d s, \d x) \right) G(x) \right] }(\theta) Q^{(2)}(\d t, \d x, \d \theta).
\end{equation*}

Let us now give the pathwise representation of the corresponding populations. From Equation (\ref{equation_Zt1_mesure}), it follows that
\begin{equation}
Z_t^{(1)}( \d a, \d x) =  \int_{(0,t]} \int_{\R_+\times \R_+} \1_{\left[0,\rho(s) H(x)\right]}(\theta)  \delta_{(t-s,x)} ( \d a, \d x)Q^{(1)}(\d s, \d x, \d \theta).
\end{equation}
and 
\begin{equation}
Z_t^{(2)}( \d a, \d x) =  \int_{(0,t]} \int_{\R_+\times \R_+} \1_{ \left[ 0,  \left( \mu(s)+ \la Z^{(2)}_{s-}, \Phi_{s} \ra +\la Z^{(1)}_{s-}, \Psi_{s} \ra \right) G(x)  \right]}  (\theta) \delta_{(t-s,x)} ( \d a, \d x)Q^{(2)}(\d s, \d x, \d \theta).
\end{equation}

Such representations are used in the field of stochastic population dynamics for populations with ages and/or characteristics (see in particular \cite{MELEARD2004}, \cite{TRAN2008} and \cite{workingpaper}). As for the standard Hawkes process, the pathwise representation has many advantages. In particular, it allows to derive the full system of equations and to identify each generation (see Section \ref{section_construction}).
More importantly, this formulation makes the link between the Hawkes process literature and the field of stochastic population dynamics. To further investigate this link seems to be a promising direction for future research.


\section*{Conclusion}
We introduced the concept of age pyramid for a class of Hawkes processes with general immigrants. The virtue of this approach is to keep track of all past events. This allows tractable computations for the Hawkes process with general immigrants whose fertility functions are time dependent generalizations of the popular exponential case, providing natural extensions of the existing results in this direction. In addition, we illustrated the pathwise construction of the Hawkes dynamics and its underlying population process.  On the whole, our approach seems to reconcile two definitions of Hawkes processes, through an intensity process or a branching dynamics. 
This framework appears to be a promising direction for further research. As an example, the large population asymptotics in the field of measure-valued population dynamics could give further insights on the macroscopic behavior of Hawkes processes.

\section*{Acknowledgements}
The author is grateful to his supervisor Nicole El Karoui for her help to improve the results and the whole paper. The author also thanks Mathieu Rosenbaum, Thibault Jaisson and Monique Jeanblanc for fruitful discussions and enlightening comments.

\section*{Appendix}

\paragraph{Proof of Corollary \ref{corollaire_laplace}.}

Let us identify the terms in Equation (\ref{equation_differentielle}). Let us denote $A_t=(A_{-1}(t),..., A_{n-1}(t))$.
 The identification of the first component gives 
 \begin{equation}
 \label{equation_A-1}
c_{-1} A_{n-1}(t) + A_{-1}^{'}(t)=0.
\end{equation} 
The second component leads to
\begin{equation}
\label{equation_B0}
A_0^{'}(t) + c_0 A_{n-1}(t) +e^{A_t. m} -1 =0.
\end{equation}
As for $1 \leq k \leq n-1$, we get
\begin{equation}
\label{equation_Bk}
A_{k-1}(t) +c_k A_{n-1}(t)+A_k^{'}(t)=0.
\end{equation}
Recursive computation of (\ref{equation_Bk}) provides for $0 \leq k \leq n-1$,
\begin{equation}
\label{equation_Bkbis}
A_k(t)=(-1)^{n-1-k} A_{n-1}^{(n-1-k)}(t) + \sum_{l=1}^{n-1-k} (-1)^l c_{k+l} A_{n-1}^{(l-1)}(t). 
\end{equation}
One deduces from (\ref{equation_Bkbis}) that
\begin{equation}
\label{equation_B0bis}
A_0'(t)=(-1)^{n-1}A_{n-1}^{(n)}(t) + \sum_{k=1}^{n-1} (-1)^k c_k A_{n-1}^{(k)}(t). 
\end{equation}
Let us intoduce the function $G(t)=\int_T^t A_{n-1}(s) \d s$ and choose $A_{-1}(t)=\theta_1 - c_{-1} G(t)$ that satisifes (\ref{equation_A-1}).
Now, put  (\ref{equation_Bkbis}) and (\ref{equation_B0bis}) into (\ref{equation_B0}) to get the following non-linear ordinary differential equation for $G$, 

\begin{equation}
\label{equation_Bnmoinsun}
\begin{split}
&(-1)^{n-1} G^{(n+1)}(t) + \sum_{k=0}^{n-1} (-1)^k c_k G^{(k+1)}(t) +\exp \bigg(\theta_1 - c_{-1} G(t) \\
&+ m_{n-1} G'(t) + \sum_{k=0}^{n-2} m_k  \left[(-1)^{n-1-k} G^{(n-k)}(t) + \sum_{l=1}^{n-1-k} (-1)^l c_{k+l} G^{(l)}(t) \right] \bigg)-1 =0.
\end{split}
\end{equation}

Let us simplify the sum in the exponential. By changing variable $k$ into $n-1-k$, it is equal to
$
 \sum_{k=1}^{n-1} m_{n-1-k}  (-1)^{k} G^{(k+1)}(t)  + \sum_{k=1}^{n-1} \sum_{l=1}^{k} (-1)^l m_{n-1-k} c_{n-1-k+l} G^{(l)}(t).
$

Then exchanging the sums leads to\\
$
  \sum_{k=1}^{n-1} m_{n-1-k}  (-1)^{k} G^{(k+1)}(t)  + \sum_{l=1}^{n-1}  (-1)^l  \left(\sum_{k=l}^{n-1}m_{n-1-k} c_{n-1-k+l} \right)  G^{(l)}(t).
$

Finally, by setting $l \leftarrow l+1$ and exchanging notations $k$ and $l$, (\ref{equation_Bnmoinsun}) becomes
\begin{equation}
\label{equation_edo_preuve}
\begin{split}
&(-1)^{n-1} G^{(n+1)}(t) + \sum_{k=0}^{n-1} (-1)^k c_k G^{(k+1)}(t)  + \exp \left( \theta_1 - c_{-1} G(t)+\sum_{k=0}^{n-1} b_k G^{(k+1)}(t)\right)-1 =0,
\end{split}
\end{equation}
where %
for $0\leq k \leq n-1$, $b_k=(-1)^k \left( m_{n-1-k} - \sum_{l=k+1}^{n-1} m_{n-1-l} c_{n-l+k} \right)$. 

Now, let us use (\ref{equation_laplace_prop}) with (\ref{equation_B0}) to get
\begin{equation*}
\begin{split}
&\E\left[ \exp \left(v.X_T \right)\right] = \exp \left(-\bmu \int_0^T (A_0^{'}(t) + c_0 A_{n-1}(t)) \d t \right),\\
&= \exp \left(-\bmu \int_0^T \left((-1)^{n-1}G^{(n+1)}(t) + \sum_{k=0}^{n-1} (-1)^k c_k G^{(k+1)}(t) \right) \d t \right),\\
&= \exp \left(-\bmu \left((-1)^{n-1}(G^{(n)}(T)- G^{(n)}(0)) + \sum_{k=0}^{n-1} (-1)^k c_k (G^{(k)}(T)-G^{(k)}(0)) \right)  \right),
\end{split}
\end{equation*}
where the second equality comes from (\ref{equation_B0bis}).
Let us set $A_0(T)=\theta_2$ and for $1\leq k \leq n-1$, $A_{k}(T)=0$. 
One can show by (\ref{equation_Bkbis}) that the previous conditions are equivalent to the terminal values $G^{(n)}(T)=(-1)^{n-1} \theta_2$ and for $1\leq k \leq n-1$, $G^{(k)}(T)=0$. Note that by definition of $G$ we also get $G(T)=0$.
We thus get 
\begin{equation*}
\begin{split}
&\E\left[ \exp \left(\theta_1 N_T + \theta_2 . \la Z_T, \phi \ra \right)\right] = \exp \left\{-\bmu \left(\theta_2 +(-1)^{n}G^{(n)}(0) + \sum_{k=0}^{n-1} (-1)^{k+1} c_k G^{(k)}(0) \right)\right\}.
\end{split}
\end{equation*}
This concludes the proof. $\diamond$

\bibliographystyle{ormsv080}
\bibliography{biblio}

\end{document}